\documentclass[12pt]{article}
\usepackage[a4paper,margin=1in]{geometry}
\usepackage{amssymb}
\usepackage{amsmath}
\usepackage{amsfonts}
\usepackage{amsthm}
\usepackage{color}
\usepackage{float}
\usepackage[all]{xy}
\usepackage{scalerel}
\usepackage{mathabx}
\usepackage{delimset}
\usepackage{mathtools}

\newcommand{\C}{\mathbb{C}}
\newcommand{\Z}{\mathbb{Z}}

\newcommand{\N}{\mathbb{N}}

\newcommand{\M}{\mathrm{M}}
\newcommand{\I}{\mathrm{I}}
\newcommand{\D}{\mathbf{D}}

\newcommand{\thp}{\Theta}

\newtheorem{theorem}{Theorem}
\newtheorem{lemma}{Lemma}
\newtheorem{definition}{Definition}
\newtheorem{proposition}{Proposition}
\newtheorem{corollary}{Corollary}
\newtheorem{example}{Example}

\newtheorem{remark}{Remark}

\begin{document}

\title{\textbf{A Resolvent Approach to Generalized Lambert Series and $q$-Series Identities}}
\author{Ronald Orozco L\'opez}

\newcommand{\Addresses}{{% additional braces for segregating \footnotesize
  \bigskip
  \footnotesize

  %R.~Orozco, \textsc{Departamento de Matematicas, Universidad de los Andes, 
  % Carrera 1 N. 18A-12 Bogot\'a, Colombia}\par\nopagebreak
  \textit{E-mail address}, R.~Orozco: \texttt{rj.orozco@uniandes.edu.co}
  
}}

\maketitle
%\tableofcontents

\begin{abstract}
We introduce a partial-theta-type \(q\)-operator
$$
\Theta(yD_q)=\sum_{n\ge0}q^{\binom n2}y^nD_q^n,
\qquad
D_qf(x)=\frac{f(qx)}{x},
$$
and show that it admits the resolvent representation
$$
\Theta(yD_q)=\left(I-\frac yxM_q\right)^{-1},
$$
where $M_qf(x)=f(qx)$. This identity provides a unified operational
framework for generalized Lambert series and their Mehler, Rogers, and
bilateral analogues. Starting from ordinary and bilateral generating
functions, we obtain Lambert-type expansions and derive consequences
involving basic hypergeometric series, Ramanujan's ${}_1\psi_1$
summation, and Kronecker-type theta identities. The method gives a
compact way to generate families of $q$-series identities from a
first-order $q$-difference resolvent.
\end{abstract}


\section{Introduction}

Lambert series occupy a classical place in the theory of \(q\)-series,
partition identities, divisor sums, and special functions
\cite{agarwal,Berndt1999Lambert,Andrews1986}. In its
standard form, a Lambert series is written as
\[
   \sum_{n\ge 1} a_n \frac{q^n}{1-q^n},
\]
where \(\{a_n\}_{n\ge1}\) is an arithmetic or combinatorial sequence.
Such series appear naturally in generating functions, in the study of
arithmetical functions, and in identities connected with the work of Ramanujan
\cite{agarwal,Berndt1991,AdigaBerndtBhargavaWatson1985}. More general Lambert-type series, involving rational
functions of \(q^n\), also arise in connection with theta functions,
basic hypergeometric series, and bilateral summations.

The purpose of this paper is to introduce a simple operational
framework for deriving generalized Lambert series and related
\(q\)-series identities. The main ingredient is a non-standard
\(q\)-difference operator
\[
   D_q f(x)=\frac{f(qx)}{x},
\]
together with the associated partial-theta-type operator
\[
   \Theta(yD_q)
   =
   \sum_{n\ge0} q^{\binom n2} y^n D_q^n.
\]
Although this operator is defined through a partial-theta-type series,
its main usefulness comes from a resolvent representation. If
\(M_q\) denotes the \(q\)-dilation operator \(M_qf(x)=f(qx)\), then
\[
   q^{\binom n2}D_q^n f(x)=x^{-n}M_q^n f(x),
\]
and hence
\[
   \Theta(yD_q)
   =
   \left(I-\frac{y}{x}M_q\right)^{-1}.
\]
Thus \(\Theta(yD_q)\) is the inverse, in a formal or convergent sense,
of the first-order \(q\)-difference operator
\[
   I-\frac{y}{x}M_q.
\]
Equivalently, it provides an operational solution of equations of the
form
\[
   F(x)-\frac{y}{x}F(qx)=G(x).
\]

This elementary observation leads to a unified treatment of several
families of Lambert-type identities. If
\[
   f(u)=\sum_{n\ge0}a_nu^n,
\]
then the action of the operator gives
\[
   \sum_{n\ge0}a_n
   \frac{x^{n+1}}{x-q^ny}z^n
   =
   \sum_{k\ge0}
   \left(\frac{y}{x}\right)^k f(q^kxz),
\]
under suitable convergence conditions, or as a formal identity whenever
the corresponding manipulations are justified algebraically. This
identity may be regarded as the basic Lambert transform associated with the operator \(\Theta(yD_q)\). Classical Lambert series are recovered by suitable choices of the parameters, while more general choices lead to Lambert-type sums with additional variables and weights.

Thus, the novelty of the approach lies not in the individual geometric expansion alone, but in the systematic use of the resolvent to pass from ordinary and bilateral generating functions to Lambert-type $q$-series identities, including Mehler, Rogers, and Kronecker-type forms. In particular, we develop two natural extensions of this construction: a Lambert–Mehler type series, involving two Lambert denominators, and a Lambert–Rogers type series, where the coefficients arise from two generating functions. These extensions show that the same resolvent mechanism produces double-sum identities and products of generating functions evaluated at $q$-dilated arguments, thereby providing a unified and compact derivation of identities that otherwise appear as separate manipulations of generating functions (see also \cite{Schmidt2020,AmdeberhanAndrewsBallantine} for related developments).

A further part of the paper is devoted to consequences involving basic
hypergeometric series. By applying the operational formulas to standard
\(q\)-series expansions, we obtain identities involving
\({}_2\phi_1\) and \({}_3\phi_2\) series. These formulas illustrate how
the Lambert-type denominators generated by the operator interact with
\(q\)-shifted factorials and basic hypergeometric summations.

The bilateral theory is especially relevant for identities of
Ramanujan type. If
$$
   f(u)=\sum_{n=-\infty}^{\infty}a_nu^n
$$
is a bilateral series convergent in an annulus, then the same operational argument gives
$$
   \sum_{n=-\infty}^{\infty}a_n
   \frac{x^{n+1}}{x-q^ny}z^n
   =
   \sum_{k\ge0}
   \left(\frac{y}{x}\right)^k f(q^kxz).
$$
This bilateral Lambert transform applies naturally to Ramanujan's
${}_1\psi_1$ summation and to theta-type coefficient sequences. In
particular, the bilateral framework gives identities involving
Ramanujan's ${}_1\psi_1$ and Kronecker-type theta identities. These applications are not independent additions; they are theta-type specializations of the same bilateral
Lambert mechanism.

The paper is organized as follows. In Section 2, we recall the basic notation for $q$-shifted factorials, basic hypergeometric series, theta functions, and Ramanujan's bilateral summation. In Section 3, we introduce the operator $D_q$, prove the formula for its iterates, and derive the resolvent representation of $\Theta(yD_q)$. Section 4 develops the generalized Lambert series associated with ordinary generating functions. Sections 5 and 6 treat the Lambert--Mehler and Lambert--Rogers extensions. Section 7 develops the bilateral Lambert transform and its theta-type applications, including Ramanujan's ${}_1\psi_1$ summation and Kronecker-type identities. Section 8 contains concluding remarks.

\section{Preliminaries}

The $q$-shifted factorial of $a$ is given by $(a;q)_{0}=1$ and
\begin{align*}
    (a;q)_{n}&=\prod_{k=0}^{n-1}(1-q^{k}a),\text{ for }n=1,2,\ldots,\\
    (a;q)_{\infty}&=\lim_{n\rightarrow\infty}(a;q)_{n}=\prod_{k=0}^{\infty}(1-aq^{k}).
\end{align*}
Following Gasper and Rahman \cite{GasperRahman}, the ${}_r\phi_{s}$ basic hypergeometric series is defined by
\begin{equation*}
    {}_r\phi_{s}\left(
    \begin{array}{c}
         a_{1},a_{2},\ldots,a_{r} \\
         b_{1},\ldots,b_{s}
    \end{array}
    ;q,z
    \right)=\sum_{n=0}^{\infty}\frac{(a_{1},a_{2},\ldots,a_{r};q)_{n}}{(q;q)_{n}(b_{1},b_{2},\ldots,b_{s};q)_{n}}\Big[(-1)^{n}q^{\binom{n}{2}}\Big]^{1+s-r}z^n
\end{equation*}
where
\begin{align*}
    (a_{1},a_{2},\ldots,a_{m};q)_{n}&=(a_{1};q)_{n}(a_{2};q)_{n}\cdots(a_{m};q)_{n}.
\end{align*}
Equivalently, 
\begin{equation*}
    (a_{1},a_{2},\ldots,a_{m};q)_{\infty}=(a_{1};q)_{\infty}(a_{2};q)_{\infty}\cdots(a_{m};q)_{\infty}.
\end{equation*}
In this paper, we will frequently use the $q$-binomial theorem
\begin{equation}\label{cauchy}
    \frac{(az;q)_{\infty}}{(z;q)_\infty}=\sum_{n=0}^{\infty}\frac{(a;q)_n}{(q;q)_n}z^n
\end{equation}
and the $q$-Gauss sum
\begin{equation}\label{eqn_qgauss}
    {}_{2}\phi_{1}\left(\begin{array}{c}
         a,b\\
         c 
    \end{array};q,\frac{c}{ab}\right)=\frac{(c/a,c/b;q)_{\infty}}{(c,c/ab;q)_{\infty}}.
\end{equation}
The following easily verified identities will be used in this paper:
\begin{align}
    (a;q)_{n}&=\frac{(a;q)_{\infty}}{(aq^n;q)_{\infty}},\label{eqn_iden1}\\
    (a;q)_{n+k}&=(a;q)_{n}(aq^{n};q)_{k},\label{eqn_iden2}\\
    (q^{-n}a;q)_{\infty}&=\frac{(-a)^n}{q^{\binom{n+1}{2}}}(q/a;q)_{n}(a;q)_{\infty}.\label{eqn2}
\end{align}
For a nonzero complex number $z$ and $q=e^{2\pi i\tau}$, where $\mathrm{Im}(\tau)>0$, the Jacobi theta function is defined as
\begin{align}
    \theta_1(z\mid \tau)&=2\sum_{n=0}^{\infty}(-1)^nq^{(2n+1)^2/8}\sin(2n+1)z.\label{theta2}
\end{align}
The Jacobi theta function has the following product representation 
\begin{align}
    \theta_1(z\mid \tau)&=2q^{1/8}(\sin z)(q;q)_\infty(qe^{2iz};q)_\infty(q/e^{2iz;q})_\infty.\label{theta_prod2}
\end{align}
We follow standard notation for Jacobi theta functions; see
\cite{WhittakerWatson,AdigaBerndtBhargavaWatson1985}.
Another famous identity is Ramanujan's summation formula
\begin{equation}\label{eqn_ramanujan}
    {}_{1}\psi_{1}(a;b;q,z)=\sum_{n=-\infty}^{\infty}\frac{(a;q)_{n}}{(b;q)_{n}}z^n=\frac{(q,b/a,az,q/az;q)_{\infty}}{(b,q/a,z,b/az;q)_{\infty}}
\end{equation}
provided that $\vert b/a\vert<\vert z\vert<1,\hspace{0.2cm}\vert q\vert<1$. For background on basic hypergeometric series and Ramanujan's \({}_1\psi_1\) summation, see \cite{GasperRahman,Andrews1986}.

Unless otherwise stated, identities involving infinite series are understood either as formal identities in the corresponding region of formal convergence or analytically under absolute convergence assumptions that justify the interchange of summation and operator action.

\section{The partial-theta-type operator and a first-order $q$-difference equation}\label{sec3}

\subsection{A non-standard $q$-derivative and the $q$-dilation operator}

\begin{definition}
We define the non-standard $q$-derivative $\mathbf{D}_{q}$ of the function $f(x)$ as
\begin{equation}
    \mathbf{D}_{q}f(x)=
    \begin{cases}
    \frac{f(q x)}{x},&\text{ if }x\neq0;\\
    \lim_{x\rightarrow0}\frac{f(q x)}{x},&\text{ if }x=0,
    \end{cases}
\end{equation}    
provided that the limit exists.
\end{definition}

\begin{theorem}\label{prop_properties}
For all $\alpha,\beta,\gamma\in\C$,
\begin{enumerate}
    \item $\mathbf{D}_{q}\{\alpha f+\beta g\}=\alpha\mathbf{D}_{q}f+\beta\mathbf{D}_{q}g$.
    \item $\mathbf{D}_{q}\{\gamma\}=\frac{\gamma}{x}$.
    \item $\mathbf{D}_{q}\{x^{n}\}=q^{n}x^{n-1}$, for $n\in\Z$.
\end{enumerate}
\end{theorem}

%\begin{remark}
    
%\end{remark}

\begin{proposition}[{\bf Product $q$-rule}]\label{prop_prod_rule}
\begin{equation}
 \mathbf{D}_{q}\{f(x)g(x)\}=f(q x)\mathbf{D}_{q}g(x)=\mathbf{D}_{q}f(x)\cdot g(q x)=x\mathbf{D}_{q}f(x)\mathbf{D}_{q}g(x).   
\end{equation}
\end{proposition}

\begin{proposition}\label{prop_dnf}
For all $n\in\N$
    \begin{equation}
        \mathbf{D}_{q}^{n}f(x)=\frac{f(q^nx)}{q^{\binom{n}{2}}x^n},
    \end{equation}
where $\D_{q}^{n}=\D_{q}\D_{q}^{n-1}$    .
\end{proposition}
\begin{proof}
The proof is by induction on $n$. For $n=1$, use the definition $\mathbf{D}_{q}$. Now suppose it is true for $n$. We will prove for $n+1$. We have
\begin{align*}
    \mathbf{D}_{q}^{n+1}f(x)=\mathbf{D}_{q}\mathbf{D}_{q}^{n}f(x)=\mathbf{D}_{q}\left\{\frac{f(q^nx)}{q^{\binom{n}{2}}x^{n}}\right\}=\frac{f(q^{n}q x)}{q^{\binom{n}{2}}q^nx^nx}=\frac{f(q^{n+1}x)}{q^{\binom{n+1}{2}}x^{n+1}}
\end{align*}
\end{proof}

\begin{proposition}[{\bf Leibniz $q$-rule}]\label{prop_leibniz}
For all $n\in\N$,
\begin{equation}
    \mathbf{D}_{q}^{n}(fg)=q^{\binom{n}{2}}x^{n}\mathbf{D}_{q}^{n}(f)\mathbf{D}_{q}^{n}(g).
\end{equation}    
\end{proposition}

\begin{proposition}
For all $n\in\N$,
\begin{equation}
    \D_{q}^n(f_{1}f_{2}\cdots f_{k})=q^{(k-1)\binom{n}{2}}x^{(k-1)n}\D_{q}^n(f_{1})\D_{q}^n(f_{2})\cdots\D_{q}^n(f_{k}).
\end{equation}
\end{proposition}

Define the $q$-dilation operator $\M_{q}$ acting on $f(x)$ to be $\M_{q}\{f(x)\}=f(q x)$. Then 
\begin{equation}
    \D_q=X^{-1}\M_q,
\end{equation}
where $X^{-1}$ is the operator of multiplication by the variable $x^{-1}$.

\begin{lemma}\label{lem_powD}
    \begin{equation}
        q^{\binom{n}{2}}\D_q^nf(x)=x^{-n}\M_q^nf(x).
    \end{equation}
\end{lemma}
\begin{proof}
From the Proposition \ref{prop_dnf}
\begin{equation}
    q^{\binom{n}{2}}\D_q^nf(x)=q^{\binom{n}{2}}\frac{f(q^nx)}{q^{\binom{n}{2}}x^n}=x^{-n}f(q^nx).
\end{equation}
\end{proof}

%\begin{remark}
    
%\end{remark}

\subsection{The partial-theta-type operator}

\begin{definition}\label{def_PTO}
We define the partial-theta-type operator associated with $\D_{q}$ by
\begin{equation*}
    \thp(y\D_{q})=\sum_{n=0}^{\infty}q^{\binom{n}{2}}y^n\D_{q}^n.
\end{equation*}
\end{definition}

\begin{proposition}\label{prop_import}
    \begin{equation}
        \thp(y\D_q)=\frac{\I}{\I-\frac{y}{x}\M_q}.
    \end{equation}
\end{proposition}
\begin{proof}
From Lemma \ref{lem_powD}
\begin{align*}
    \thp(y\D_{q})f(x)=\sum_{k=0}^{\infty}q^{\binom{k}{2}}y^{k}\D_{q}^kf(x)=\sum_{k=0}^{\infty}\left(\frac{y}{x}\right)^k\M_q^kf(x)=\frac{\I}{\I-\frac{y}{x}\M_{q}}f(x).
\end{align*}
\end{proof}

Lemma 1 shows that $\thp(y\D_q)$ reduces to a geometric series in $\M_q$. 

\subsection{Operational inversion and $q$-difference equations}

A fundamental property of the partial-theta-type operator $\thp(y\D_q)$ is its role as the inverse of the first-order $q$-difference operator $\I-\frac{y}{x}\M_q$. Specifically, we have 
\begin{equation}
    (\I-\frac{y}{x}\M_q)\thp(y\D_q)=\I,
\end{equation}
which implies that $\thp(y\D_q)$ provides a formal solution to equations of the form
\begin{equation}
    f(x)-\frac{y}{x}f(qx)=g(x).
\end{equation}
This operational inversion allows us to systematically construct explicit solutions for a wide class of $q$-difference equations. In particular, the operator $\thp(y\D_q)$ generates series representations with denominators of the Lambert type, thus unifying the treatment of classical and generalized Lambert series within an operational framework.

\begin{example}[\bf Monomial right-hand side]\label{exa_lambert}
Consider the equation
\begin{equation}
    f(x)-\frac{y}{x}f(qx)=x^n
\end{equation}
for $n\in\N$. Applying the operator $\thp(y\D_q)$, we obtain the explicit solution
\begin{align}
        f(x)=\thp(y\D_{q})\big\{x^{n}\big\}&=\frac{x^{n+1}}{x-q^ny}.
    \end{align}
This formula shows that the operator produces rational functions with denominators of the Lambert type, which are central to the theory of Lambert series. In general, the solution of the $q$-difference equation
\begin{equation}
    f(x)-\frac{y}{x}f(qx)=x^{\alpha n+\beta}
\end{equation}
for all $\alpha,\beta\in\N$, is
    \begin{equation}
        f(x)=\frac{x^{\alpha n+\beta+1}}{x-q^{\alpha n+\beta}y}.
    \end{equation}
\end{example}

\begin{example}[\bf General power series]
Let $g(x)=\sum_{n=0}^{\infty}a_nx^n$ be an analytic function. Then, the solution to
\begin{equation}
    f(x)-\frac{y}{x}f(qx)=g(x)
\end{equation}
is given by
\begin{equation}
    f(x)=\thp(y\D_q)\{g(x)\}=\sum_{n=0}^{\infty}a_n\frac{x^{n+1}}{x-q^ny}.
\end{equation}
Alternatively, using the operational formula for $\thp(y\D_q)$, we have
\begin{equation}
    f(x)=\sum_{k=0}^{\infty}\left(\frac{y}{x}\right)^kg(q^kx).
\end{equation}
This dual representation links the operator approach to classical generating-function manipulations.
\end{example}

\begin{example}[\bf Classical Lambert series]
As a special case, let $g(x)=\sum_{n=1}^{\infty}a_nx^n$ and set $y=0$. Then,
\begin{equation}
    f(x)=\sum_{n=1}^{\infty}a_nx^n,
\end{equation}
which recovers the original generating function. However, for $y=x=q$, we obtain
\begin{equation}
    f(x)=\sum_{n=1}^{\infty}a_n\frac{q^{n+1}}{q-q^nq}=\sum_{n=1}^{\infty}a_n\frac{q^{n}}{1-q^n},
\end{equation}
which is the classical Lambert series. This demonstrates how the operator framework naturally generalizes and recovers classical results.
\end{example}

\section{Generalized Lambert series}

\begin{definition}
Let $a=\{a_n\}$ be a sequence. For $x\neq y$, we define the generalized Lambert series related to the sequence $a$ as
    \begin{equation}
        \mathcal{L}_a(x,y,z;q)=\sum_{n=0}^{\infty}a_{n}\frac{x^{n+1}}{x-q^ny}z^n    .
    \end{equation}
\end{definition}

\begin{remark}
The following elementary specializations recover several classical Lambert-type denominators:
\begin{itemize}
    \item 
    \begin{equation}
        \mathcal{L}_{a}(1,\pm q,q;q)=\mathcal{L}_{a}(q,\pm q^2,1;q)=\sum_{n=0}^{\infty}a_{n}\frac{q^n}{1\mp q^{n+1}}.
    \end{equation}
\item 
\begin{equation}
        \mathcal{L}_{a}(q,\pm q^2,q;q)=\sum_{n=0}^{\infty}a_{n}\frac{q^{2n}}{1\mp q^{n+1}}.
    \end{equation}
\item     
\begin{equation}
        \mathcal{L}_{a}(1,\pm q^{-\beta},q^\alpha;q^\alpha)=\sum_{n=0}^{\infty}a_{n}\frac{q^{\alpha n}}{1\mp q^{\alpha n+\beta}}.
    \end{equation}    
\end{itemize}
\end{remark}

\begin{theorem}\label{theo_gls}
Let \(f(u)=\sum_{n\ge0}a_nu^n\) be a power series with radius of convergence \(R\). Then
\begin{equation}
    \mathcal{L}_{a}(x,y,z;q)=\sum_{k=0}^{\infty}\left(\frac{y}{x}\right)^kf(q^kxz),
\end{equation}
provided that \(|y|<|x|\), \(|xz|<R\), and \(0<|q|<1\).
\end{theorem}
\begin{proof}
From Proposition \ref{prop_import}
    \begin{align*}
        \mathcal{L}_{a}(x,y,z;q)&=\sum_{n=0}^{\infty}a_{n}\frac{x^{n+1}}{x-q^ny}z^n\\
        &=\sum_{n=0}^{\infty}a_{n}\thp(y\D_{q})\{x^n\}z^n
        =\thp(y\D_{q})\left\{\sum_{n=0}^{\infty}a_{n}(xz)^n\right\}
        =\sum_{k=0}^{\infty}\left(\frac{y}{x}\right)^kf(q^kxz).
    \end{align*}
\end{proof}

\begin{corollary}[\bf Cauchy formula and Lambert series]\label{cor_CGL}
    \begin{equation}
        \sum_{n=0}^{\infty}\frac{(a;q)_n}{(q;q)_n}\frac{x^{n+1}}{x-q^ny}z^n=\frac{(axz;q)_\infty}{(xz;q)_\infty}{}_{2}\phi_{1}\left(
    \begin{array}{c}
         xz,q \\
         axz
    \end{array}
    ;q,\frac{y}{x}
    \right).
    \end{equation}
\end{corollary}
\begin{proof}
By using the Theorem \ref{theo_gls},
\begin{align*}
    \sum_{n=0}^{\infty}\frac{(a;q)_n}{(q;q)_n}\frac{x^{n+1}}{x-q^ny}z^n&=\sum_{k=0}^{\infty}\left(\frac{y}{x}\right)^k\frac{(q^kaxz;q)_{\infty}}{(q^kxz;q)_\infty}\\
    &=\frac{(axz;q)_{\infty}}{(xz;q)_\infty}\sum_{k=0}^{\infty}\left(\frac{y}{x}\right)^k\frac{(xz;q)_{k}}{(axz;q)_k}\\
    &=\frac{(axz;q)_\infty}{(xz;q)_\infty}{}_{2}\phi_{1}\left(
    \begin{array}{c}
         xz,q \\
         axz
    \end{array}
    ;q,\frac{y}{x}
    \right).
\end{align*}
\end{proof}

\section{Generalized Lambert-Mehler Series}

\begin{definition}
The generalized Lambert-Mehler type series is defined by
\begin{equation}
    \mathcal{L}_{a}(x,y;z,w;t;q)=\sum_{n=0}^{\infty}a_{n}\frac{x^{n+1}}{x-q^ny}\cdot\frac{z^{n+1}}{z-q^nw}t^n,
\end{equation}
where $x\neq y$ and $z\neq w$.
\end{definition}

\begin{remark}
The following elementary specializations recover several classical Lambert-type denominators:
\begin{itemize}
    \item 
    \begin{equation}
    \mathcal{L}_{a}(1,\pm y^{-\beta};1,\pm q^{-\gamma};q^\alpha;q^{\alpha})=\sum_{n=0}^{\infty}a_{n}\frac{q^{2\alpha n}}{(1\mp q^{\alpha n-\beta})(1\mp q^{\alpha n-\gamma})}.
\end{equation}
If $\beta=\gamma$, then
\begin{equation}
    \mathcal{L}_{a}(1,\pm y^{-\beta};1,\pm q^{-\beta};q^\alpha;q^{\alpha})=\sum_{n=0}^{\infty}a_{n}\frac{q^{2\alpha n}}{(1\mp q^{\alpha n-\beta})^2}.
\end{equation}
\item 
\begin{equation}
    \mathcal{L}_{a}(1,\pm y^{-\beta};1,\mp q^{-\gamma};q^\alpha;q^{\alpha})=\sum_{n=0}^{\infty}a_{n}\frac{q^{2\alpha n}}{(1\mp q^{\alpha n-\beta})(1\pm q^{\alpha n-\gamma})}.
\end{equation}
If $\beta=\gamma$, then
\begin{equation}
    \mathcal{L}_{a}(1,\pm y^{-\beta};1,\mp q^{-\beta};q^\alpha;q^{\alpha})=\sum_{n=0}^{\infty}a_{n}\frac{q^{2\alpha n}}{1- q^{2\alpha n-2\beta}}.
\end{equation}
\end{itemize}
\end{remark}

\begin{theorem}[{\bf Lambert-Mehler formula}]\label{mehler}
Let \(f(u)=\sum_{n\ge0}a_nu^n\) be a power series with radius of convergence \(R\). Then 
    \begin{equation}
        \mathcal{L}_{a}(x,y;z,w;t;q)=\sum_{k=0}^{\infty}\sum_{l=0}^{\infty}\left(\frac{w}{z}\right)^k\left(\frac{y}{x}\right)^lf(q^{k+l}txz)
    \end{equation}
provided that $0<\vert y\vert<\vert x\vert$, $0<\vert w\vert<\vert z\vert$, $\vert txz\vert<R$, and $0<\vert q\vert<1$.    
\end{theorem}
\begin{proof}
    \begin{align*}
        \sum_{n=0}^{\infty}a_{n}\frac{x^{n+1}}{x-q^ny}\cdot\frac{z^{n+1}}{z-q^nw}t^n&=\sum_{n=0}^{\infty}a_{n}\thp(y\D_{q})\{x^n\}\frac{z^{n+1}}{z-q^nw}t^n\\
        &=\thp(y\D_{q})\left\{\sum_{n=0}^{\infty}a_{n}\frac{z^{n+1}}{z-q^nw}(tx)^n\right\}\\
        &=\thp(y\D_{q})\left\{\sum_{k=0}^{\infty}\left(\frac{w}{z}\right)^kf(q^ktxz)\right\}\\
        &=\sum_{k=0}^{\infty}\left(\frac{w}{z}\right)^k\thp(y\D_{q})\{f(q^ktxz)\}\\
        &=\sum_{k=0}^{\infty}\left(\frac{w}{z}\right)^k\sum_{l=0}^{\infty}\left(\frac{y}{x}\right)^lf(q^{k+l}txz).
    \end{align*}
\end{proof}
Let $\thp^{(x)}(y\D_q)$ denote the operator $\thp(y\D_q)$ acting on variable $x$. Then
\begin{equation}
    \mathcal{L}_{a}(x,y;z,w;t;q)=\thp^{(x)}(y\D_q)\thp^{(z)}(w\D_q)f(txz).
\end{equation}

\begin{corollary}[\bf Cauchy formula and Lambert-Mehler series]\label{cor_CGML}
    \begin{multline}
        \sum_{n=0}^{\infty}\frac{(a;q)_n}{(q;q)_n}\frac{x^{n+1}}{x-q^ny}\cdot\frac{z^{n+1}}{z-q^nw}t^n\\=\frac{(atxz;q)_\infty}{(txz;q)_\infty}\sum_{k=0}^{\infty}\left(\frac{w}{z}\right)^k\frac{(txz;q)_k}{(atxz;q)_k}{}_{2}\phi_{1}\left(
    \begin{array}{c}
         q^ktxz,q \\
         q^katxz
    \end{array}
    ;q,\frac{y}{x}
    \right)
    \end{multline}
\end{corollary}
\begin{proof}
We use the Theorem \ref{mehler}
    \begin{align*}
        &\sum_{k=0}^{\infty}\sum_{l=0}^{\infty}\left(\frac{w}{z}\right)^k\left(\frac{y}{x}\right)^l\frac{(q^{k+l}atxz;q)_\infty}{(q^{k+l}txz;q)_\infty}\\
        &\hspace{1cm}=\frac{(atxz;q)_\infty}{(txz;q)_\infty}\sum_{k=0}^{\infty}\sum_{l=0}^{\infty}\left(\frac{w}{z}\right)^k\left(\frac{y}{x}\right)^l\frac{(txz;q)_{k+l}}{(atxz;q)_{k+l}}.
     \end{align*}
Now, by using the identity $(a;q)_{k+l}=(a;q)_k(q^ka;q)_l$     
     \begin{align*}
        &\sum_{k=0}^{\infty}\sum_{l=0}^{\infty}\left(\frac{w}{z}\right)^k\left(\frac{y}{x}\right)^l\frac{(q^{k+l}atxz;q)_\infty}{(q^{k+l}txz;q)_\infty}\\
        &\hspace{1cm}=\frac{(atxz;q)_\infty}{(txz;q)_\infty}\sum_{k=0}^{\infty}\left(\frac{w}{z}\right)^k\frac{(txz;q)_k}{(atxz;q)_k}\sum_{l=0}^{\infty}\left(\frac{y}{x}\right)^l\frac{(q^ktxz;q)_{l}}{(q^katxz;q)_{l}}\\
        &\hspace{1cm}=\frac{(atxz;q)_\infty}{(txz;q)_\infty}\sum_{k=0}^{\infty}\left(\frac{w}{z}\right)^k\frac{(txz;q)_k}{(atxz;q)_k}{}_{2}\phi_{1}\left(
    \begin{array}{c}
         q^ktxz,q \\
         q^katxz
    \end{array}
    ;q,\frac{y}{x}
    \right).
    \end{align*}
\end{proof}

\section{Generalized Lambert-Rogers series}

\begin{definition}
The generalized Lambert-Rogers type series is defined as
\begin{equation}
    \mathcal{L}_{a,b}(x,y,t,s;q)=\sum_{n=0}^{\infty}\sum_{m=0}^{\infty}a_{n}b_{m}\frac{x^{n+m+1}}{x-q^{n+m}y}t^ns^m,
\end{equation}
for $x\neq y$.
\end{definition}
As a specialization, we have
\begin{equation}
    \mathcal{L}_{a,b}(1,\pm y^{-\beta},q^\alpha,q^\alpha;q^\alpha)=\sum_{n=0}^{\infty}\sum_{m=0}^{\infty}a_{n}b_{m}\frac{q^{\alpha(n+m)}}{1\mp q^{\alpha(n+m)-\beta}}.
\end{equation}

\begin{theorem}[{\bf Lambert-Rogers formula}]\label{rogers}
Let $f(u)=\sum_{n\ge0}a_nu^n$ and $g(u)=\sum_{n\ge0}b_nu^n$ be power series with radius of convergence $R_1$ and $R_2$, respectively. Then
    \begin{equation}
        \mathcal{L}_{a,b}(x,y,t,s;q)=\sum_{k=0}^{\infty}\left(\frac{y}{x}\right)^{k}f(q^ktx)g(q^ksx)
    \end{equation}
provided that $0<\vert y\vert<\vert x\vert$ and $\vert xt\vert<R_{1}$, $\vert sx\vert<R_{2}$, and $0<\vert q\vert<1$.    
\end{theorem}
\begin{proof}
    \begin{align*}
        \sum_{n=0}^{\infty}\sum_{m=0}^{\infty}a_{n}b_{m}\frac{x^{n+m+1}}{x-q^{n+m}y}t^ns^m
        &=\sum_{n=0}^{\infty}\sum_{m=0}^{\infty}a_{n}b_{m}\thp(y\D_{q})\{x^{n+m}\}t^ns^m\\
        &=\thp(y\D_{q})\left\{\sum_{n=0}^{\infty}\sum_{m=0}^{\infty}a_{n}b_{m}(xt)^n(xs)^m\right\}\\
        &=\thp(y\D_{q})\left\{f(tx)g(sx)\right\}\\
        &=\sum_{k=0}^{\infty}\left(\frac{y}{x}\right)^{k}f(q^ktx)g(q^ksx).
    \end{align*}
\end{proof}
From Theorem \ref{rogers},
\begin{equation}
    \mathcal{L}_{a,b}(x,y,t,s;q)=\thp(y\D_q)\{f(tx)g(sx)\}.
\end{equation}

\begin{corollary}[\bf Cauchy formula and Lambert-Rogers series]\label{cor_CGRL}
    \begin{equation}
        \sum_{n=0}^{\infty}\sum_{m=0}^{\infty}\frac{(a;q)_n}{(q;q)_n}\frac{(b;q)_m}{(q;q)_m}\frac{x^{n+m+1}}{x-q^{n+m}y}t^ns^m=\frac{(atx,bsx;q)_\infty}{(tx,sx;q)_\infty}
    {}_{3}\phi_{2}\left(
    \begin{array}{c}
         tx,sx,q \\
         atx,bsx
    \end{array}
    ;q,\frac{y}{x}
    \right).
    \end{equation}
\end{corollary}
\begin{proof}
From Theorem \ref{rogers},
\begin{align*}
    &\sum_{n=0}^{\infty}\sum_{m=0}^{\infty}\frac{(a;q)_n}{(q;q)_n}\frac{(b;q)_m}{(q;q)_m}\frac{x^{n+m+1}}{x-q^{n+m}y}t^ns^m\\
    &\hspace{1cm}=\sum_{k=0}^{\infty}\left(\frac{y}{x}\right)^{k}\frac{(aq^ktx;q)_\infty}{(q^ktx;q)_\infty}\frac{(bq^ksx;q)_\infty}{(q^ksx;q)_\infty}\\
    &\hspace{1cm}=\frac{(atx;q)_\infty}{(tx;q)_\infty}\frac{(bsx;q)_\infty}{(sx;q)_\infty}\sum_{k=0}^{\infty}\left(\frac{y}{x}\right)^{k}\frac{(tx;q)_k}{(atx;q)_k}\frac{(sx;q)_k}{(bsx;q)_k}\\
    &\hspace{1cm}=\frac{(atx,bsx;q)_\infty}{(tx,sx;q)_\infty}
    {}_{3}\phi_{2}\left(
    \begin{array}{c}
         tx,sx,q \\
         atx,bsx
    \end{array}
    ;q,\frac{y}{x}
    \right).
\end{align*}
\end{proof}

\section{Bilateral Lambert Series and Theta-Type Applications}

The bilateral Lambert transform introduced in this section is the natural
counterpart of the ordinary Lambert transform in Sections~4--6.
Its main advantage is that it converts \emph{theta-type bilateral generating functions} into Lambert denominators \(1-q^n(\cdot)\), producing Ramanujan-type identities. In particular, Ramanujan's ${}_1\psi_1$ summation appears as a canonical bilateral
input and Kronecker-type identities arise as theta-type specializations of the same mechanism.

We first state the bilateral Lambert series (Theorem~\ref{theo_BGL}).
We then apply it to Ramanujan's \({}_1\psi_1\) kernel (Corollary~\ref{coro_R-L}),
extend it to two-denominator and product settings (Theorems~\ref{bimehler}--\ref{birogers}),
and finally specialize to theta-type coefficients leading to Kronecker-type theta identities.

\subsection{Bilateral generalized Lambert series and Ramanujan’s ${}_{1}\Psi_{1}$ summation }

\begin{definition}[{\bf Bilateral generalized Lambert series}]
Let $\{a_n\}_{n\in\mathbb{Z}}$ be a doubly infinite sequence.
We define the bilateral generalized Lambert series associated with $\{a_n\}$ by
\begin{equation}
\sum_{n=-\infty}^{\infty}
a_n\frac{x^{n+1}}{x-q^ny}z^n,
\end{equation}
with $x\neq y$.
\end{definition}

\begin{theorem}[\bf Bilateral Lambert representation]\label{theo_BGL}
Let $f(u) = \sum_{n=-\infty}^{\infty}a_nu^n$ be a bilateral series convergent in an annulus. Then
\begin{equation}
    \sum_{n=-\infty}^{\infty}a_n \frac{x^{n+1}}{x-q^ny}z^n=
    \sum_{k=0}^{\infty}\left(\frac{y}{x}\right)^kf(q^k xz).
\end{equation}
\end{theorem}
\begin{proof}
\begin{align*}
    \sum_{n=-\infty}^{\infty}a_n \frac{x^{n+1}}{x-q^ny}z^n&=\thp(y\D_q)\left\{f(zx)\right\}=\sum_{k=0}^{\infty}\left(\frac{y}{x}\right)^kf(q^k xz).
\end{align*}   
\end{proof}

\begin{lemma}\label{lemma2}
For all $n\in\N$, and $x\neq0$, $a\neq0$, $b\neq0$, we have
    \begin{equation}
        \M_{q}^n\left\{\frac{(ax,q/ax;q)_{\infty}}{(x,b/ax;q)_{\infty}}\right\}=q^nx^n\frac{(x;q)_{n}}{(qax/b;q)_{n}}\frac{(ax,q/ax;q)_{\infty}}{(x,b/ax;q)_{\infty}}
    \end{equation}
provided that $\max\{\vert x\vert,\vert b/ax\vert\}<1$.    
\end{lemma}
\begin{proof}
From Proposition \ref{prop_dnf}, and if we apply Eqs. (\ref{eqn_iden1}) and (\ref{eqn2}), we obtain
    \begin{align*}
        \M_{q}^n\left\{\frac{(ax,q/ax;q)_{\infty}}{(x,b/ax;q)_{\infty}}\right\}&=\frac{(q^nax,q^{1-n}/ax;q)_{\infty}}{(q^nx,q^{-n}b/ax;q)_{\infty}}\\
        &=\frac{(x;q)_{n}(-1)^n(ax;q)_{n}(ax)^nq^{\binom{n+1}{2}}}{a^nq^{\binom{n}{2}}(ax;q)_{n}(-1)^n(qax/b;q)_{n}}\frac{(ax,q/ax;q)_{\infty}}{(x,b/ax;q)_{\infty}}\\
        &=q^nx^n\frac{(x;q)_{n}}{(qax/b;q)_{n}}\frac{(ax,q/ax;q)_{\infty}}{(x,b/ax;q)_{\infty}}.
    \end{align*}
This proves the lemma.
\end{proof}

In the remainder of this section, we apply Theorem~\ref{theo_BGL} with
bilateral inputs \(f\) of hypergeometric or theta type.

\begin{corollary}[\bf Ramanujan summation formula and bilateral Lambert series]\label{coro_R-L}
\begin{equation}
    \sum_{n=-\infty}^{\infty}\frac{(a;q)_n}{(b;q)_n}\frac{x^{n+1}}{x-q^ny}z^n={}_{1}\Psi_{1}(a;b;q,zx){}_{2}\phi_{1}\left(
    \begin{array}{c}
         xz,q \\
         qazx/b
    \end{array}
    ;q,qy
    \right).
\end{equation}
\end{corollary}
\begin{proof}
From Theorem \ref{theo_BGL}
    \begin{align*}
        \sum_{n=-\infty}^{\infty}\frac{(a;q)_n}{(b;q)_n}\frac{x^{n+1}}{x-q^ny}z^n
        =\frac{(q,b/a;q)_{\infty}}{(b,q/a;q)_{\infty}}\sum_{n=0}^{\infty}\left(\frac{y}{x}\right)^n\M_q^n\left\{\frac{(azx,q/azx;q)_{\infty}}{(zx,b/azx;q)_{\infty}}\right\}.
    \end{align*}
From Lemma \ref{lemma2},
    \begin{align*}
        \sum_{n=-\infty}^{\infty}\frac{(a;q)_n}{(b;q)_n}\frac{x^{n+1}}{x-q^ny}z^n
        &=\frac{(q,b/a;q)_{\infty}}{(b,q/a;q)_{\infty}}\sum_{n=0}^{\infty}\left(\frac{y}{x}\right)^nq^nx^n\frac{(xz;q)_{n}}{(qaxz/b;q)_{n}}\frac{(axz,q/axz;q)_{\infty}}{(xz,b/axz;q)_{\infty}}\\
        &={}_{1}\Psi_{1}(a;b;q,zx)\sum_{n=0}^{\infty}\frac{(xz;q)_{n}}{(qaxz/b;q)_{n}}(qy)^n\\
        &={}_{1}\Psi_{1}(a;b;q,zx){}_{2}\phi_{1}\left(
    \begin{array}{c}
         xz,q \\
         qazx/b
    \end{array}
    ;q,qy
    \right).
    \end{align*}
\end{proof}

\begin{remark}[\bf Product specialization]\label{rem:productSpecialization}
If $qy=a/b$ in Corollary~\ref{coro_R-L}, then one obtains the explicit product
\begin{equation}
\sum_{n=-\infty}^{\infty}\frac{(a;q)_n}{(b;q)_n}\,
\frac{bx^{n+1}}{bx-q^{n-1} a}\,z^n
=
{}_{1}\Psi_{1}(a;b;q,xz)\frac{b-azx}{b-a},
\end{equation}
where the final simplification uses the \(q\)-Gauss sum.
\end{remark}

\medskip
\noindent\textbf{Summary of the mechanism.}
Theorem~\ref{theo_BGL} provides a single template:
\[
\text{(bilateral input \(f\))} \quad \Longrightarrow \quad
\text{Lambert denominators \(1-q^n(\cdot)\)}.
\]
The next subsections show how this template extends to two denominators
(Mehler-type) and to products (Rogers-type), and then yields theta-type applications.

\subsection{Bilateral Mehler and Rogers extensions}

\begin{definition}
The bilateral generalized Lambert-Mehler type series is defined by
\begin{equation}
    \sum_{n=-\infty}^{\infty}a_{n}\frac{x^{n+1}}{x-q^ny}\cdot\frac{z^{n+1}}{z-q^nw}t^n,
\end{equation}
where $x\neq y$ and $z\neq w$.
\end{definition}

\begin{theorem}[{\bf Bilateral Lambert-Mehler formula}]\label{bimehler}
Let $f(u) = \sum_{n=-\infty}^{\infty}a_nu^n$ be a bilateral series convergent in an annulus. Then
    \begin{equation}
        \sum_{n=-\infty}^{\infty}a_{n}\frac{x^{n+1}}{x-q^ny}\cdot\frac{z^{n+1}}{z-q^nw}t^n=\sum_{k=0}^{\infty}\sum_{l=0}^{\infty}\left(\frac{w}{z}\right)^k\left(\frac{y}{x}\right)^lf(q^{k+l}txz)
    \end{equation}
\end{theorem}

\begin{corollary}[\bf Ramanujan ${}_{1}\Psi_{1}$ and bilateral Lambert-Mehler series]\label{coroBLM}
For parameters in the region of convergence,
    \begin{multline}
        \sum_{n=-\infty}^{\infty}\frac{(a;q)_n}{(b;q)_n}\frac{x^{n+1}}{x-q^ny}\cdot\frac{z^{n+1}}{z-q^nw}t^n\\
        ={}_{1}\Psi_{1}(a;b;q,txz)\sum_{k=0}^{\infty}\left(qwtx\right)^k\frac{(txz;q)_k}{(qatxz/b;q)_k}{}_{2}\phi_{1}\left(
    \begin{array}{c}
         q^ktxz,q \\
         q^{k+1}atxz/b
    \end{array}
    ;q,qtyz
    \right).
    \end{multline}
\end{corollary}
\begin{proof}
From Theorem \ref{bimehler}
    \begin{align*}
        &\sum_{n=-\infty}^{\infty}\frac{(a;q)_n}{(b;q)_n}\frac{x^{n+1}}{x-q^ny}\cdot\frac{z^{n+1}}{z-q^nw}t^n\\
        &\hspace{1cm}=\sum_{k=0}^{\infty}\sum_{l=0}^{\infty}\left(\frac{w}{z}\right)^k\left(\frac{y}{x}\right)^l\frac{(q,b/a,aq^{k+l}txz,q^{1-k-l}/atxz;q)_{\infty}}{(b,q/a,q^{k+l}txz,q^{-k-l}b/atxz;q)_{\infty}}\\
    \end{align*}
Now, from Lemma \ref{lemma2}
\begin{align*}
    &\sum_{n=-\infty}^{\infty}\frac{(a;q)_n}{(b;q)_n}\frac{x^{n+1}}{x-q^ny}\cdot\frac{z^{n+1}}{z-q^nw}t^n\\
    &\hspace{1cm}=\frac{(q,b/a;q)_\infty}{(b,q/a;q)_\infty}\sum_{k=0}^{\infty}\sum_{l=0}^{\infty}\left(\frac{w}{z}\right)^k\left(\frac{y}{x}\right)^l(qtxz)^{k+l}\frac{(txz;q)_{k+l}}{(qatxz/b;q)_{k+l}}\frac{(atxz,q/atxz;q)_{\infty}}{(txz,b/atxz;q)_{\infty}}\\
        &\hspace{1cm}={}_{1}\Psi_{1}(a;b;q,txz)\sum_{k=0}^{\infty}\sum_{l=0}^{\infty}\left(\frac{w}{z}\right)^k\left(\frac{y}{x}\right)^l(qtxz)^{k+l}\frac{(txz;q)_{k+l}}{(qatxz/b;q)_{k+l}}\\
        &\hspace{1cm}={}_{1}\Psi_{1}(a;b;q,txz)\sum_{k=0}^{\infty}\left(qwtx\right)^k\frac{(txz;q)_k}{(qatxz/b;q)_k}\sum_{l=0}^{\infty}\left(qytz\right)^l\frac{(q^ktxz;q)_{l}}{(q^{k+1}atxz/b;q)_{l}}\\
        &\hspace{1cm}={}_{1}\Psi_{1}(a;b;q,txz)\sum_{k=0}^{\infty}\left(qwtx\right)^k\frac{(txz;q)_k}{(qatxz/b;q)_k}{}_{2}\phi_{1}\left(
    \begin{array}{c}
         q^ktxz,q \\
         q^{k+1}atxz/b
    \end{array}
    ;q,qtyz
    \right).
\end{align*}
\end{proof}

\begin{definition}
The generalized bilateral Lambert-Rogers type series is defined as
\begin{equation}
    \sum_{n=-\infty}^{\infty}\sum_{m=-\infty}^{\infty}a_{n}b_{m}\frac{x^{n+m+1}}{x-q^{n+m}y}t^ns^m,
\end{equation}
for $x\neq y$.
\end{definition}

\begin{theorem}[{\bf Bilateral Lambert-Rogers formula}]\label{birogers}
Let $f(u)=\sum_{n=-\infty}^{\infty}a_nu^n$ and $g(u)=\sum_{n=-\infty}^{\infty}b_nu^n$ be bilateral series convergent in an annulus. Then
    \begin{equation}
        \sum_{n=-\infty}^{\infty}\sum_{m=-\infty}^{\infty}a_{n}b_{m}\frac{x^{n+m+1}}{x-q^{n+m}y}t^ns^m=\sum_{k=0}^{\infty}\left(\frac{y}{x}\right)^{k}f(q^ktx)g(q^ksx).
    \end{equation}
\end{theorem}

\begin{corollary}[\bf Ramanujan ${}_{1}\Psi_{1}$ and bilateral Lambert-Rogers series]\label{coroBLR}
For parameters in the region of convergence,
\begin{multline}
    \sum_{n=-\infty}^{\infty}\sum_{m=-\infty}^{\infty}\frac{(a;q)_n}{(b;q)_n}\frac{(c;q)_m}{(d;q)_m}\frac{x^{n+m+1}}{x-q^{n+m}y}t^ns^m\\
    ={}_{1}\Psi_{1}(a;b;q,tx){}_{1}\Psi_{1}(c;d;q,sx){}_{3}\phi_{2}\left(
    \begin{array}{c}
         tx,sx,q \\
         qatx/b,qcsx/d
    \end{array}
    ;q,q^2tsxy
    \right).
\end{multline}
\end{corollary}
\begin{proof}
From Theorem \ref{birogers}
    \begin{align*}
        &\sum_{n=-\infty}^{\infty}\sum_{m=-\infty}^{\infty}\frac{(a;q)_n}{(b;q)_n}\frac{(c;q)_m}{(d;q)_m}\frac{x^{n+m+1}}{x-q^{n+m}y}t^ns^m\\
        &=\frac{(q,b/a;q)_{\infty}}{(b,q/a;q)_{\infty}}\frac{(q,d/c;q)_{\infty}}{(d,q/c;q)_{\infty}}\sum_{n=0}^{\infty}\left(\frac{y}{x}\right)^n\M_q^n\left\{\frac{(atx,q/atx;q)_{\infty}}{(tx,b/atx;q)_{\infty}}\frac{(csx,q/csx;q)_{\infty}}{(sx,d/csx;q)_{\infty}}\right\}.
    \end{align*}
Now, we use Lemma \ref{lemma2}    
\begin{align*}
    &\sum_{n=-\infty}^{\infty}\sum_{m=-\infty}^{\infty}\frac{(a;q)_n}{(b;q)_n}\frac{(c;q)_m}{(d;q)_m}\frac{x^{n+m+1}}{x-q^{n+m}y}t^ns^m=\frac{(q,b/a;q)_{\infty}}{(b,q/a;q)_{\infty}}\frac{(q,d/c;q)_{\infty}}{(d,q/c;q)_{\infty}}\\
    &\hspace{1cm}\times\sum_{n=0}^{\infty}\left(\frac{y}{x}\right)^n\frac{q^n(tx)^n(tx;q)_{n}(atx,q/atx;q)_{\infty}}{(qatx/b;q)_{n}(tx,b/atx;q)_{\infty}}\frac{q^n(sx)^n(sx;q)_{n}(csx,q/csx;q)_{\infty}}{(qcsx/d;q)_{n}(sx,d/csx;q)_{\infty}}\\
    &\hspace{1cm}={}_{1}\Psi_{1}(a;b;q,tx){}_{1}\Psi_{1}(c;d;q,sx)\sum_{n=0}^{\infty}\left(q^2tsxy\right)^n\frac{(tx;q)_{n}}{(qatx/b;q)_{n}}\frac{(sx;q)_{n}}{(qcsx/d;q)_{n}}\\
    &\hspace{1cm}={}_{1}\Psi_{1}(a;b;q,tx){}_{1}\Psi_{1}(c;d;q,sx){}_{3}\phi_{2}\left(
    \begin{array}{c}
         tx,sx,q \\
         qatx/b,qcsx/d
    \end{array}
    ;q,q^2tsxy
    \right).
\end{align*}  
\end{proof}

\subsection{Kronecker-type identities}

We next indicate how Kronecker-type theta identities arise from the same bilateral framework. The key point is that Kronecker identities express quotients of theta functions as bilateral sums with Lambert denominators \(1-q^n(\cdot)\), matching the
kernel in Theorem~\ref{theo_BGL}.

\medskip
\noindent\textbf{Step 1: a \({}_1\psi_1\) specialization.}
Taking \(b=qa\) in Ramanujan's \({}_1\psi_1\) summation and using
\((a;q)_n/(aq;q)_n=(1-a)/(1-aq^n)\), we obtain, for \(a\neq0\) not an integral power
of \(q\) and \(|q|<|z|<1\),
\begin{equation}\label{eq:KroneckerKernel}
\sum_{n=-\infty}^{\infty}\frac{z^n}{1-aq^n}
=
\frac{(q;q)_\infty^2\,(az,q/az;q)_\infty}{(a,q/a,z,q/z;q)_\infty}.
\end{equation}

\medskip
\noindent\textbf{Step 2: conversion to theta quotients.}
Replacing \(z\mapsto e^{2iy}\) and \(a\mapsto e^{2iz}\) in \eqref{eq:KroneckerKernel},
and using the infinite product representation of \(\theta(z\mid \tau)\) from the
preliminaries, we arrive at Kronecker's classical identity. The Kronecker theta function is naturally connected with Ramanujan's
\({}_1\psi_1\) summation and theta quotients; see, for example,
\cite{Kronecker1881,Liu2021Kronecker}.
\begin{equation}\label{eq:KroneckerIdentity}
\sum_{n=-\infty}^{\infty}\frac{e^{2niy}}{1-q^n e^{2iz}}
=
\frac{i}{2}\,K_y(z\mid\tau),
\end{equation}
where the Kronecker function is
\begin{equation}\label{eq:KroneckerFunction}
K_y(z\mid\tau)
=
\frac{\theta_1^{\prime}(0\mid\tau)\,\theta_1(z+y\mid\tau)}
{\theta_1(z\mid\tau)\,\theta_1(y\mid\tau)},
\qquad
\theta_1(y\mid\tau)\neq0,\;\theta_1(z\mid\tau)\neq0.
\end{equation}

\medskip
\noindent\textbf{Step 3: extensions via the bilateral Mehler/Rogers mechanism.}
Using Corollary \ref{coro_R-L} as the basic Kronecker kernel and applying the extensions in Corollaries \ref{coroBLM} and \ref{coroBLR}, we obtain:
\begin{align}
\sum_{n=-\infty}^{\infty}\frac{e^{2niy}}{(1-q^n e^{2iz})(1-q^{n}\alpha)}
&=
\frac{i}{2(1-q\alpha)}\,K_y(z\mid\tau),
\\[2mm]
\sum_{n=-\infty}^{\infty}\frac{e^{2niy}}{(1-q^n e^{2iz})(1-q^n\alpha)(1-q^n\beta)}
&=
\frac{iK_y(z\mid\tau)}{2(1-q\alpha e^{2iy})(1-q\beta e^{2iy})},
\\[2mm]
\sum_{n=-\infty}^{\infty}\sum_{m=-\infty}^{\infty}
\frac{e^{2i(nz+mw)}}{(1-q^n e^{2iz})(1-q^m e^{2ix})(1-q^{n+m}\alpha)}
&=
-\frac{K_y(z\mid\tau)\,K_x(w\mid\tau)}{4(1-q^2\alpha e^{2i(z+w)})}.
\end{align}

\medskip
\noindent
These formulas illustrate how theta quotients (Kronecker-type objects) emerge as
theta-type specializations of the bilateral Lambert framework.

\section{Concluding remarks}

We have shown that the partial-theta-type operator $\Theta(yD_q)$, associated with the non-standard $q$-difference operator $D_qf(x)=f(qx)/x$, admits a simple resolvent representation in terms of the $q$-dilation operator. This representation gives a unified mechanism for generating generalized Lambert series, Lambert--Mehler and Lambert--Rogers extensions, and bilateral theta-type identities. The applications to Ramanujan's ${}_1\psi_1$ and Kronecker-type identities illustrate that the same operator acts as a bridge between ordinary generating functions, bilateral summations, and theta quotients.

Possible further directions include applying the same resolvent method to other bilateral summations, elliptic analogues, or Lambert-type series arising from modular and mock modular objects.

\end{document}